\documentclass[11pt]{article}
\usepackage{setspace}
\usepackage{amsmath}
\usepackage{amsthm}
\usepackage{amssymb}
\usepackage{extarrows}
\usepackage{graphicx}
\usepackage{sidecap}
\usepackage[T1]{fontenc}
\usepackage{microtype}
\usepackage{floatrow}
\usepackage[english]{babel}
\usepackage[hmarginratio=1:1,top=32mm,columnsep=20pt]{geometry}
\usepackage[hang, small,labelfont=bf,up,textfont=it,up]{caption}
\usepackage{booktabs}
\usepackage{indentfirst}
\usepackage{lettrine}
\usepackage{subfigure}
\usepackage{diagbox}
\usepackage{multirow}
\usepackage{multicol}
\usepackage{enumitem}
\setlist[itemize]{noitemsep}
\usepackage{abstract}

\usepackage{titlesec}
\usepackage[ruled]{algorithm2e}
\usepackage{titling}
\usepackage{hyperref}
\usepackage{tkz-euclide,pgfplots}

\newtheorem{theorem}{Theorem}

\newtheorem{assumption}{Assumption}

\newtheorem{remark}{Remark}

\pretitle{\begin{center}\LARGE\bfseries}
\posttitle{\end{center}}

\title{Robust parameter estimation of regression models under weakened moment assumptions}

\author{Kangqiang Li\thanks{Corresponding author E-mail address: 11935023@zju.edu.cn (Kangqiang Li)} \qquad Songqiao Tang\thanks{E-mail address: 11835013@zju.edu.cn (Songqiao Tang)}\qquad Lixin Zhang\thanks{E-mail address: stazlx@zju.edu.cn (Lixin Zhang)}\\[1ex] 
\normalsize School of Mathematical Sciences, Zhejiang University,  Hangzhou,  Zhejiang 310027,  China \\ 
}
\date{}

\begin{document}
\maketitle
\begin{abstract}
This paper provides some extended results on estimating parameter matrix of several regression models when the covariate or response possesses weaker moment condition. We study the $M$-estimator of Fan et al. (Ann Stat 49(3):1239--1266, 2021) for matrix completion model with $(1+\epsilon)$-th moment noise. The corresponding phase transition phenomenon is observed. When $1> \epsilon>0$, the robust estimator possesses a slower convergence rate compared with previous literature. For high dimensional multiple index coefficient model, we propose an improved estimator via applying the element-wise truncation method to handle heavy-tailed data with finite fourth moment. The extensive simulation study validates our theoretical results.
\end{abstract}

\section{\textbf{Introduction}}
Under the traditional settings, sub-Gaussian assumption is often required for noise and design in regression problems. Due to the heavy-tailed phenomena of real-world data, in recent years, there has been a growing body of literature on the robust regression estimation when the covariate and response are heavy-tailed. For linear-type models, Fan et al. (2017)\cite{symmetry} applied the Huber (1964)\cite{location parameter}'s loss to the sparse linear model and showed that under $(2+\epsilon)$-th moment assumption on the noise, the proposed estimator exhibits the same statistical error rate as that of the light-tail noise case. Further, Sun et al. (2020)\cite{Adaptive} proposed the adaptive Huber regression and extended the result of Fan et al. (2017)\cite{symmetry} to the case of $(1+\epsilon)$-th moment condition on the noise. A tight phase transition for the estimation error of the regression parameter was established which parallelled those first discovered by Bubeck et al. (2013)\cite{Bandits with heavy tail} and Devroye et al. (2016)\cite{Devroye} for robust mean estimation without finite variance. Motivated by Sun et al. (2020)\cite{Adaptive}, Tan et al. (2022)\cite{reduced} established the similar phase transition results on sparse reduced rank regression. Fan et al. (2021)\cite{shrinkage} focused on robust estimation of the trace regression and their $M$-estimator achieves the minimax statistical error rate under only bounded $(2+\epsilon)$-th moment response or both bounded fourth moment design and response. Afterwards, Han et al (2021)\cite{post-selection} constructed a post-selection inference procedure via the Huber loss for high-dimensional linear model and Zhang (2021)\cite{series} investigated Huber robust estimators for high-dimensional heavy-tailed time series. Avella-Medina et al. (2018)\cite{precision matrices} and Ke et al (2019)\cite{User-friendly} applied Huber loss to construct robust covariance and precision matrix estimators. Zhu and Zhou (2020)\cite{Taming} studied the corrupted general linear model with heavy-tailed data under finite fourth moment assumption.\par
For robust parameter estimation of sparse non-linear regression problem, Neykov et al. (2016)\cite{Neykov} analyzed least squares with $L_1$ penalization in high-dimensional single index model (SIM) under Gaussian designs. Plan and Vershynin (2016)\cite{Plan2016} and Plan et al. (2017)\cite{Plan2017} considered high-dimensional SIM under Gaussian and elliptical designs. Yang et al. (2018)\cite{thresholded} proposed a robust estimator of high-dimensional SIM when the covariate and response only have the bounded fourth moment. The proposed estimator achieves the optimal error bound via the truncation procedure. Furthermore, Goldstein et al. (2018)\cite{Structured} analyzed high-dimensional SIM under heavy-tailed elliptical distribution. Na and Kolar (2021)\cite{volatility} developed an estimation procedure of the parametric components in high-dimensional index volatility models under finite moments condition. Fan et al. (2022)\cite{Understanding} studied implicit regularization in SIM with general heavy-tailed data. As an extension of SIM, Na et al. (2019)\cite{varying} considered high dimensional varying index coefficient model introduced by Ma and Song (2015)\cite{Varying}. For estimating sparse parameter matrix, they required the existence of bounded $6$-th moment of the design and response in order to obtain the optimal rate, but whether the moment constraint can be further relaxed is unknown.  Meanwhile, it is worth noting that Fan et al. (2021)\cite{shrinkage} tested the superiority of their estimator for trace regression via selecting the scaled Cauchy noise, beyond the corresponding theoretical condition. Motivated by those, a natural question arises:\par
\emph{Can we further generalize their results and obtain the sharp estimation rates?}\par
To address this problem, on the basis of Fan et al. (2021)\cite{shrinkage}'s work, we further study matrix completion model in which the noise distribution has no finite variance. The applicable condition of their $M$-estimator is broadened. Simultaneously, the smooth phase transition of the convergence rate is also observed. As a generalization of matrix completion model, we also consider robust parameter estimation of high-dimensional varying index coefficient model. To handle heavy-tailed data with only finite fourth moment, we give a robust element-wise truncated estimator (see (\ref{eq2.5})) based on the research of Na et al. (2019)\cite{varying}. It turns out that under finite fourth moment assumption, our method shows the robustness against the low order moments. The proposed estimator can achieve the same statistical error rate as that of Na et al. (2019)\cite{varying} with finite fourth moment. Meanwhile, the data-driven method facilitates to calibrate truncation parameters and is more convenient than the cross-validation method of Na et al. (2019)\cite{varying}. Note that this paper investigates the mean estimation with homoscedastic noise, since bounded moments condition has been allowed for the variance estimation.\par
The remainder of our paper is organized as follows. In Section \ref{sec3}, we analyze two specific regression problems and derive the statistical error rates of the corresponding $M$-estimators under weaker moment assumptions. In Section \ref{sec4}, some numerical simulations on synthetic data are presented and show an agreement with the theoretical results. Concluding remarks are drawn in Section \ref{sec5}. All the proofs are presented in the supplementary material.
\subsection*{\textbf{Notation}}\label{SecDef}
For any positive integer $n$, we denote the set $\{1,2,\ldots,n\}$ by $[n]$. For two matrices $X, Y \in \mathbb{R}^{d_{1}\times d_{2}}$, $\langle X,Y\rangle:=\text{tr}(X^{T}Y)$. For a matrix $A=(a_{ij}) \in \mathbb{R}^{d_{1}\times d_{2}}$, the max norm and Frobenius norm of $A$ are defined as $\|A\|_{\max}=\max_{i \in [d_{1}],j \in [d_{2}]}|a_{i,j}|$ and $\|A\|_{F}=\sqrt{\sum_{i\in [d_{1}],j\in [d_{2}]}a_{i,j}^2}$ respectively. $\|A\|_{\star}=\text{tr}\left(\sqrt{A^TA}\right)$, $\|A\|_{\text{op}}=\sqrt{\lambda_{\max}\left(A^TA\right)}$, $\|A\|_{1,1}=\sum_{i\in[d_{1}]}
\sum_{j\in [d_{2}]}|a_{i,j}|$, $\|A\|_{\infty}=\max_{i\in[d_{1}]}\sum_{j\in [d_{2}]}|a_{i,j}|$ and $\|A\|_{L_{1}}=\max_{j\in[d_{2}]}\sum_{i\in [d_{1}]}|a_{i,j}|$. Given two positive sequences $\{a_{n}\}_{n=1}^{\infty}$ and $\{b_{n}\}_{n=1}^{\infty}$, we use the notation $a_{n} \asymp b_{n}$, if $b_{n} \lesssim a_{n}\lesssim b_{n}$ where $a_{n}\lesssim b_{n}$ means that there exists a positive constant $C$ such that $a_{n}\leq C b_{n}$ for all $n$.

\section{\textbf{Parameter matrix estimation of linear and nonlinear statistical models}}\label{sec3}
In this section, we analyze two types of regression models and present statistical rates of the corresponding regularized estimators under weakened moment assumptions.
\subsection{\textbf{Low-rank matrix completion model with weaker moment}}\label{sec3.1}
We first consider the following matrix completion model:
\begin{equation}\label{eq4.0}
y=\left\langle X,\Theta^{\star}\right\rangle+\varepsilon
\end{equation}
where $X$ is uniformly sampled from $\{\sqrt{d_{1}d_{2}}\cdot e_{j}e_{k}^T\}_{j\in[d_{1}],
k\in[d_{2}]}$ and $\mathbb{E}(\varepsilon|X)=0$. To recover the parameter matrix $\Theta^{\star}$ under near low-rank assumption, Fan et al. (2021)\cite{shrinkage} studied the following $M$-estimator of $\Theta^{\star}$:
$$
\widehat{\Theta} = \underset{\|\Theta\|_{\max}\leq R/\sqrt{d_{1}d_{2}}}{\operatorname{argmin}}\left\{\operatorname{vec}(\Theta)^{T} \widehat{\Sigma}_{XX} \operatorname{vec}(\Theta)-2\left\langle\widehat{\Sigma}_{y X}, \Theta\right\rangle+\lambda\|\Theta\|_{\star}\right\}\emph{}
$$
where $\widehat{\Sigma}_{XX}=\frac{1}{n} \sum_{i=1}^{n} \operatorname{vec}(X_{i}) \operatorname{vec}(X_{i})^{T}$ and $\widehat{\Sigma}_{yX}=\frac{1}{n} \sum_{i=1}^{n} \operatorname{sign}\left(y_{i}\right)\left(\left|y_{i}\right| \wedge \tau\right) X_{i}$ with a truncation parameter $\tau$. Under finite $(2+\epsilon)$-th moment condition on the response, their robust estimator has the theoretical statistical error rate of order $\sqrt{r(d_{1}\vee d_{2})\log(d_{1}+d_{2})/n}$ under Frobenius norm, which is the same as that of Negahban and Wainwright (2012)\cite{convexity} for sub-exponential noise. The following theorem further relaxes the distributional conditions from the bounded $(2+\epsilon)$-th moment to $(1+\epsilon)$-th moment assumption to fill the gap for the robust estimator's scope of use.\par

\begin{theorem}\label{theorem1}
Suppose the following conditions hold:\par
(1) $\left\|\Theta^{*}\right\|_{F}\leq 1$, $\left\|\Theta^{*}\right\|_{\max}\leq R/
\sqrt{d_{1}d_{2}}$, $\left\|\Theta^{*}\right\|_{\max } /\left\|\Theta^{*}\right\|_{F}
 \leq R / \sqrt{d_{1} d_{2}}$ and $\operatorname{rank}\left(\Theta^{\star}\right)\leq r$;\par
(2) $\{X_{i}\}_{i=1}^{n}$ are i.i.d. uniformly sampled from $\{\sqrt{d_{1}d_{2}}\cdot e_{j}e_{k}^T\}_{j\in[d_{1}], k\in[d_{2}]}$ and $\mathbb{E}\left(\mathbb{E}(|\varepsilon_{i}
|^{\alpha}\big{|}X_{i})^{\log(d_{1}\vee d_{2})}\right) \leq M_{\alpha}<\infty$ for some $\alpha \in (1,2]$.\par Then for any $\delta>1$, choose $\tau\asymp \left(\frac{L_{\alpha}n}
{(d_{1}\vee d_{2})\log(d_{1}+d_{2})}\right)^{\frac{1}{\alpha}}$ and for some constant $C>0$,
$$
\lambda=4C\left(\frac{(d_{1}\vee d_{2})\log(d_{1}+d_{2})}
{n}\right)^{\frac{\alpha-1}{\alpha}}\left(L_{\alpha}^{\frac{1}{\alpha}}
\delta+R\delta+L_{\alpha}^{\frac{1}{\alpha}}\right),
$$
there exist constants $\{C_{i}\}_{i=1}^{4}$ such that as long as $n\geq(d_{1}\vee d_{2})\log(d_{1}+d_{2})$, we have with the probability at least $1-2(d_{1}+d_{2})^{1-\delta}
-C_{1}\exp\left(-C_{2}(d_{1}+d_{2})\right)$,
$$
\begin{aligned}
&\left\|\widehat{\Theta}-\Theta^{\star}\right\|_{F} \leq C_{3} \sqrt{r} \left(\frac{(d_{1}\vee d_{2})\log(d_{1}+d_{2})}{n}\right)^{\frac{\alpha-1}
{\alpha}}\left(L_{\alpha}^{\frac{1}{\alpha}}\delta
+R\delta+L_{\alpha}^{\frac{1}{\alpha}}\right), \\ &\left\| \widehat{\Theta}
-\Theta^{*}\right\|_{\star} \leq C_{4} r \left(\frac{(d_{1}\vee d_{2})\log(d_{1}+
d_{2})}{n}\right)^{\frac{\alpha-1}{\alpha}}\left(L_{\alpha}^{\frac{1}{\alpha}}
\delta+R\delta+L_{\alpha}^{\frac{1}{\alpha}}\right)
\end{aligned}
$$
where $L_{\alpha}=2^{\alpha-1}\left(R^{\alpha}+e{M_{\alpha}^{1/\log(d_{1}\vee d_{2})}}\right)$.
\end{theorem}

\begin{remark}\label{remark1}
According to Theorem \ref{theorem1}, we obtain that for some $\alpha>1$, choosing $\tau\asymp \min\{L_{\alpha}^{\frac{1}{\alpha}},
L_{2}^{\frac{1}{2}}\}\left(\frac{n}{(d_{1}\vee d_{2})\log(d_{1}+d_{2})}
\right)^{\max\{\frac{1}{\alpha},\frac{1}{2}\}}$ and $\lambda\asymp\min\{L_{\alpha}^{\frac{1}{\alpha}},L_{2}^{\frac{1}{2}}\}\left(\frac{(d_{1}\vee d_{2})\log(d_{1}+d_{2})}{n}\right)^{\min\{\frac{\alpha-1}{\alpha},\frac{1}{2}\}}$,
$$
\left\|\widehat{\Theta}-\Theta^{\star}\right\|_{F}\lesssim \sqrt{r}\min\{L_{\alpha}^{\frac{1}{\alpha}},L_{2}^{\frac{1}{2}}\} \left(\frac{(d_{1}\vee d_{2})\log(d_{1}+d_{2})}{n}\right)^{\min\{\frac{\alpha-1}
{\alpha},\frac{1}{2}\}} \;\; \text{with high probability.}
$$
Compared to the result of Fan et al. (2021)\cite{shrinkage}, when $\alpha<2$, there exists a smooth phase transition phenomenon for the statistical error rate which is in line with linear regression in Sun et al. (2020)\cite{Adaptive} and mean estimation in Bubeck et al. (2013)\cite{Bandits with heavy tail}. The truncation parameter $\tau$ and the regularized parameter $\lambda$ adapt to the moment of the noise. However, this transition is observed in the low-rank matrix completion
model via the shrinkage technique, which is a visible difference with previous literature.
\end{remark}
\begin{remark}\label{remark2}
If $\text{vec}(X)$ is a sub-Gaussian vector, the phase transition phenomenon still holds for matrix compressed sensing and multitask regression of Fan et al. (2021)\cite{shrinkage}. Specifically, when $\mathbb{E}\left(\mathbb{E}(|\varepsilon
|^{\alpha}\big{|}X)^k\right) \leq M_{\alpha}$ for some $\alpha \in (1,2]$, $k>1$ and $d_{1}+d_{2}\leq n$, by choosing $\tau\asymp \left(\frac{M_{\alpha}^{1/k}n}{d_{1}+d_{2}}\right)^{\frac{1}{\alpha}}$ and $\lambda \asymp M_{\alpha}^{\frac{1}{k\alpha}} \left(\frac{d_{1}+d_{2}}{n}\right)^{\frac{\alpha-1}
{\alpha}}$, we have that $\left\|\widehat{\Theta}-\Theta^{\star}\right\|_{F}
\lesssim \sqrt{r} M_{\alpha}^{\frac{1}{k\alpha}}\left(\frac{d_{1}+d_{2}}{n}
\right)^{\frac{\alpha-1}{\alpha}}$ and $\left\|\widehat{\Theta}-\Theta^{\star}
\right\|_{\star}\lesssim rM_{\alpha}^{\frac{1}{k\alpha}}\left(\frac{d_{1}
+d_{2}}{n}\right)^{\frac{\alpha-1}{\alpha}}$ with high probability towards matrix compressed sensing. Our simulation study confirms the above inference and the proof is omitted for less redundancy.
\end{remark}

\subsection{\textbf{High-dimensional varying index coefficient model}}\label{sec3.2}
As a generalization of model (\ref{eq4.0}), in this subsection, we concentrate on robustly estimating the direction of parameters estimation of the following varying index coefficient model:
$$
y=\sum_{i=1}^{d_{2}} z_{i} \cdot f_{i}\left(\left\langle X,\theta_{i}^{\star}\right\rangle\right)+\varepsilon
$$
where $X\in \mathbb{R}^{d_{1}}$ and $Z=(z_{1}, z_{2}, \ldots, z_{d_{2}})^{T}\in \mathbb{R}^{d_{2}}$ are independent covariates, and $\varepsilon$ is the stochastic error with $\mathbb{E}[\varepsilon\mid X,Z]=0$. We assume that $\|\theta_{i}^{\star}\|_{2}=1$ for model identifiability and $X$ has the known probability density function $p(X)$.\par
Further, assume that the following two conditions hold:\par
\begin{assumption}\label{assumption1}
Assume that the covariate $X$ has the differentiable density function $p(X): \mathbb{R}^{d_{1}} \rightarrow \mathbb{R}$ and the link functions $\{f_{i}(\cdot)| i \in [d_{2}]\} $ are differentiable such that $\mu_{i}^{\star}:=\mathbb{E}[f_{i}^{\prime}\left(\left\langle X, \theta_{i}^{\star}\right\rangle\right)]\ne 0$ for $\forall i \in [d_{2}]$, and $\mathbb{E}[Z]=0_{d_{2}\times 1}$. Denote $\Sigma^{\star}:=\mathbb{E}\left[Z Z^{T}\right]$ and $\Omega^{\star}:=\left({\Sigma}^{\star}\right)^{-1}.$ We assume ${\Omega}^{\star} \in \left\{\Omega: \Omega \succ 0,\|\Omega\|_{L_{1}}\leq \varpi, \max _{1 \leq i \leq d_2} \sum_{j=1}^{d_2}\left|(\Omega)_{i,j}\right|^{q} \leq K\right\}$ for some $\varpi$ and $q\in [0,1)$.
\end{assumption}

\begin{assumption}\label{assumption2}
There exists an absolute constant $M>0$ such that
$$\mathbb{E}\left[y^{4}\right] \vee \mathbb{E}\left[[S(X)]_{j}^{4}\right] \vee \mathbb{E}\left[z_{k}^{4}\right] \leq M, \quad \forall j \in\left[d_{1}\right], k \in\left[d_{2}\right]$$
where the first-order score function $S: \mathbb{R}^{d_{1}} \rightarrow \mathbb{R}^{d_{1}}$ is defined as $S(X):=-\nabla p(X) / p(X)$.
\end{assumption}
Based on the above assumptions and first-order Stein's identity (Stein et al. (2004)\cite{exchangeable}), according to Na et al. (2019)\cite{varying}, we have
$$
\mathbb{E}\left[y \cdot S(X) Z^{T}\right] \Omega^{\star} =\sum_{j=1}^{d_{2}}
\mathbb{E}\left[f_{j}\left(\left\langle\theta_{j}^{\star},X\right\rangle\right)
S(X)\right] \mathbb{E}\left[z_{j}\cdot Z^{T}\right] \Omega^{\star}=\sum_{j=1}^{d_{2}} \mu_{j}^{\star}\theta_{j}^{\star}e_{j}^{T}:=\left(\widetilde{\theta}_{1}, \ldots, \widetilde{\theta}_{d_{2}}\right)=\widetilde{\Theta}.
$$\par
Therefore, a feasible method to estimate the direction of $\left\{\theta_{i}^{\star}\right\}_{i=1}^{d_{1}}$ without the knowledge of the link functions $\{f_{i}(\cdot)\}_{i\in [d_{2}]}$ is pointed out. Given $n$ i.i.d. samples $\{y_{i}, X_{i}, Z_{i}\}_{i=1}^{n}$, Na et al. (2019)\cite{varying} proposed to separately truncate the data $\{y_{i}, S(X_{i}), Z_{i}\}_{i=1}^{n}$ via the function $\check{x}=x1_{\{|x|\leq \tau\}}$. By assuming 6-th moments of the covariate and response exist, their robust estimator achieved the convergence rate of order $\sqrt{sd_{2}} \left(\frac{\log (d_{1}d_{2})}{n}\right)^{\frac{1}{2}}$ under Frobenius norm. In order to further relax the moment condition, we consider $y_{i}S(X_{i})Z_{i}^{T}$ as a matrix-valued data and then use $\psi_{\tau}(x)=(|x| \wedge \tau) \operatorname{sign}(x)$ to truncate each entry of the matrix-variate data. Specifically, the robust element-wise truncated matrix estimator is defined as\par
\begin{equation}\label{eq2.5}
\widehat{\Theta}=\underset{\Theta \in \mathbb{R}^{d_{1}\times d_{2}}}
{\operatorname{argmin}}\left\{\|\Theta\|_{F}^{2}-2\left\langle\frac{1}{n}\sum_{i=1}^{n}
\psi_{\Gamma_{1}}\left(y_{i}S(X_{i})Z_{i}^{T}\right)
\widehat{\Omega},\Theta\right\rangle+\lambda\|\Theta\|_{1,1}\right\}
\end{equation}
where $\Gamma_{1}=\left(\tau_{j,k}^{(1)}\right)_{j\in[d_{1}]}^{k\in[d_{2}]}$ is a truncation parameter matrix and $\widehat{\Omega}$ is obtained by Cai et al. (2011)\cite{constrained}'s CLIME procedure:
\begin{equation}\label{eq1.6}
\widehat{\Omega}=\text{argmin} \|\Omega\|_{1,1}\;\; \text{s.t.}\;\;\left\|\widehat{\Sigma}_{n} \Omega-I_{d_{2}}\right\|_{\max }\leq\gamma,
\end{equation}
where $(\widehat{\Sigma}_{n})_{j,k}=\frac{1}{n}\sum_{i=1}^{n}\psi_{\tau_{j,k}^{(2)}}
\left(z_{j}^{(i)}z_{k}^{(i)}\right)$.

The following theorem gives the statistical error rate of the robust estimator above.
\begin{theorem}\label{theorem2}
Suppose Assumption \ref{assumption1} and \ref{assumption2} hold with $\|\theta_{j}^{\star}\|_{0} = s$ for all $j \in\left[d_{2}\right].$ For $j\in [d_{1}]$ and $k,s \in [d_{2}]$, choose $\tau_{j,k}^{(1)}\asymp M^{\frac{3}{4}}\sqrt{\frac{n}{\log (d_{1}d_{2})}}$, $\tau_{k, s}^{(2)}\asymp M^{\frac{1}{2}}
\sqrt{\frac{n}{\log d_{2}}}$, $\gamma\asymp M^{\frac{1}{2}}\varpi\sqrt{\frac{\log d_{2}}{n}}$ and
$$
\lambda = 8M^{\frac{3}{4}}\left\|\Omega^{\star}\right\|_{1,1} \sqrt{\frac{3\log (d_{1}d_{2})}{n}}+16\max _{j \in\left[d_{2}\right]}
|\mu_{j}^{\star}| \cdot\left\|\Theta^{\star}\Sigma^{\star}\right\|_{\infty} M^{\frac{1}{2}}\varpi^2\sqrt{\frac{4\log d_{2}}{n}}.
$$
Then with the probability at least $1-\frac{2}{(d_{1}d_{2})^2}-\frac{1}{d_{2}^2}-\frac{1}{d_{2}^3}$, we have
$$
\left\|\widehat{\Theta}-\widetilde{\Theta}\right\|_{F}\leq 2\lambda\sqrt{sd_{2}} \;\; \text{and} \;\; \left\|\widehat{\Theta}-\widetilde{\Theta}\right\|_{1,1}\leq 8\lambda sd_{2}.
$$
\end{theorem}

\begin{remark}\label{remark3}
From the above result, we have with high probability,
$$
\left\|\widehat{\Theta}-\widetilde{\Theta}\right\|_{F}\lesssim
\sqrt{sd_{2}} \left(\frac{\log (d_{1}d_{2})}{n}\right)^{\frac{1}{2}}\;\;\text{and} \;\; \left\|\widehat{\Theta}-\widetilde{\Theta}\right\|_{1,1}\lesssim
sd_{2} \left(\frac{\log (d_{1}d_{2})}{n}\right)^{\frac{1}{2}}
$$
which shows that the proposed estimator possesses the same statistical error rate as that of Na et al. (2019)\cite{varying} with bounded 6-th moment assumption.
\end{remark}

\section{\textbf{Simulation Study}}\label{sec4}
In this section, we provide some numerical experiments to confirm the statistical error rates of the estimators established in previous section.\par
\begin{figure}[H]
\centering
\subfigure{
\begin{tikzpicture}
\pgfplotsset{every axis x label/.append style={at={(0.5,0.1)}}, every axis y label/.append style={at={(0.15,0.5)}}}
\begin{axis}[title={$t_{2}/5$ noise}, xlabel = $\log(n)$,ylabel = {$\log(\text{error})$},title style={font=\tiny},font=\tiny,height=5.6cm, width=5.55cm]
\addplot[
    color=blue,
    mark=triangle,
    densely dashed
    ]
    coordinates {
(7.600902,-0.7668609)(8.29405,-1.019195)(8.699515,-1.283023)(8.987197,-1.396262)
(9.210340,-1.457021)(9.392662,-1.484671)(9.615805,-1.671114) };
\addplot[
   color=blue,
    mark=triangle,
    ]
    coordinates {
(7.600902,-0.9035341)(8.29405,-1.337842)(8.699515,-1.561083)(8.987197,-1.716486)
(9.210340,-1.835491)(9.392662,-1.928022)(9.615805,-2.048171) };
\addplot[
    color=red,
    mark=square,
    densely dashed
    ]
    coordinates {
(7.600902,-0.4113717)(8.29405,-0.7401078)(8.699515,-0.8970367)(8.987197,-0.9730487)
(9.210340,-1.0649986)(9.392662,-1.1055611)(9.615805,-1.1971407) };
\addplot[
    color=red,
    mark=square
    ]
    coordinates {
(7.600902,-0.4530596)(8.29405,-0.8843236)(8.699515,-1.113059)(8.987197,-1.277079)
(9.210340,-1.401099)(9.392662,-1.484733)(9.615805,-1.599713)};
\addplot[
    color=green,
    mark=x,
    densely dashed
    ]
    coordinates {
(7.600902,-0.1905578)(8.29405,-0.5272047)(8.699515,-0.7181409)(8.987197,-0.8354264)
(9.210340,-0.9070177)(9.392662,-0.9992199)(9.615805,-1.0491984) };
\addplot[
    color=green,
    mark=x
    ]
    coordinates {
(7.600902,-0.2135059)(8.29405,-0.6016985)(8.699515,-0.8656561)(8.987197,-1.0299410)
(9.210340,-1.1487168)(9.392662,-1.2362806)(9.615805,-1.3276064)
 };
\end{axis}
\end{tikzpicture}}
\subfigure{
\begin{tikzpicture}
\pgfplotsset{every axis x label/.append style={at={(0.5,0.1)}}, every axis y label/.append style={at={(0.15,0.5)}}}
\begin{axis}[title={$t_{1.5}/10$ noise},xlabel = $\log(n)$,ylabel = {$\log(\text{error})$},font=\tiny,height=5.6cm, width=5.55cm]
\addplot[
    color=blue,
    mark=triangle,
    densely dashed
    ]
    coordinates {
(8.29405,-0.7962847)(8.699515,-0.9169874)(8.987197,-0.9964296)(9.210340,-1.0238205)
(9.392662,-1.0859727)(9.615805,-1.1337300)(9.90348,-1.184304)};
\addplot[
    color=red,
    mark=square,
    densely dashed
    ]
    coordinates {
(8.29405,-0.5387244)(8.699515,-0.6205695)(8.987197,-0.6546069)(9.210340,-0.6601381)
(9.392662,-0.6096973)(9.615805,-0.6271368)(9.90348,-0.6619467)};
\addplot[
    color=green,
    mark=x,
    densely dashed
    ]
    coordinates {
(8.29405,-0.4406177)(8.699515,-0.5581774)(8.987197,-0.6210987)(9.210340,-0.5404835)
(9.392662,-0.5795378)(9.615805,-0.5453645)(9.90348,-0.480379)};
\addplot[
    color=blue,
    mark=triangle
    ]
    coordinates {
(8.29405,-1.411550)(8.699515,-1.543718)(8.987197,-1.640522)(9.210340,-1.726966)
(9.392662,-1.803851)(9.615805,-1.882606)(9.90348,-1.979254)};
\addplot[
    color=red,
    mark=square
    ]
    coordinates {
(8.29405,-1.030356)(8.699515,-1.2408279)(8.987197,-1.3284904)(9.210340,-1.4067099)
(9.392662,-1.4672619)(9.615805,-1.5432922)(9.90348,-1.646918)
};
\addplot[
    color=green,
    mark=x
    ]
    coordinates{
(8.29405,-0.7208968)(8.699515,-1.0009297)(8.987197,-1.1407467)(9.210340,-1.2501543)
(9.392662,-1.3191536)(9.615805,-1.3908870)(9.90348,-1.480585)};
\end{axis}
\end{tikzpicture}}
\subfigure{
\begin{tikzpicture}
\pgfplotsset{every axis x label/.append style={at={(0.5,0.1)}}, every axis y label/.append style={at={(0.22,0.5)}}}
\begin{axis}[title={$t_{1.1}/15$ noise},xlabel = $\log(n)$,ylabel = {$\log(\text{error})$},font=\tiny,height=5.6cm, width=5.55cm]
\addplot[
    color=blue,
    mark=triangle,
    densely dashed
    ]
    coordinates {(8.29405,0.4803439)(8.699515,0.5385629)(8.987197,0.5244956)
(9.210340,0.5557851)(9.392662,0.4685363)(9.615805,0.6197065)(9.90348,0.5266995)};
\addplot[
    color=red,
    mark=square,
    densely dashed
    ]
    coordinates {(8.29405,0.4984608)(8.699515,0.7191086)(8.987197,0.8876047)
(9.210340,0.9173844)(9.392662,1.0283437)(9.615805,0.9863225)(9.90348,1.037773)};
\addplot[
    color=green,
    mark=x,
    densely dashed
    ]
    coordinates {(8.29405,0.3850764)(8.699515,0.5624511)(8.987197,0.7581182)
(9.210340,0.7546030)(9.392662,0.8555738)(9.615805,0.9926937)(9.90348,1.174252)};
\addplot[
    color=blue,
    mark=triangle
    ]
    coordinates {
(8.29405,-1.2856682)(8.699515,-1.3537285)(8.987197,-1.3227117)(9.210340,-1.3258087)
(9.392662,-1.3507907)(9.615805,-1.3693051)(9.90348,-1.388292)};
\addplot[
    color=red,
    mark=square
    ]
    coordinates {
(8.29405,-0.846641)(8.699515,-1.1323715)(8.987197,-1.2091790)(9.210340,-1.2425272)
(9.392662,-1.2529427)(9.615805,-1.2565747)(9.90348,-1.24089)};
\addplot[
    color=green,
    mark=x
    ]
    coordinates{
(8.29405,-0.6151333)(8.699515,-0.8673865)(8.987197,-1.0167960)(9.210340,-1.1435276)
(9.392662,-1.2072011)(9.615805,-1.2370208)(9.90348,-1.237318)};
\end{axis}
\end{tikzpicture}}
\includegraphics[width=0.95\columnwidth]{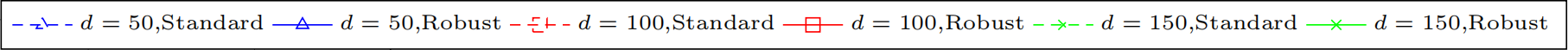}
\caption{Matrix completion: The $x$-axis and $y$-axis represent logarithmic sample size and $\log\left\|\widehat{\Theta}-\Theta^{*}\right\|_{F}$.}\label{fig1}
\end{figure}
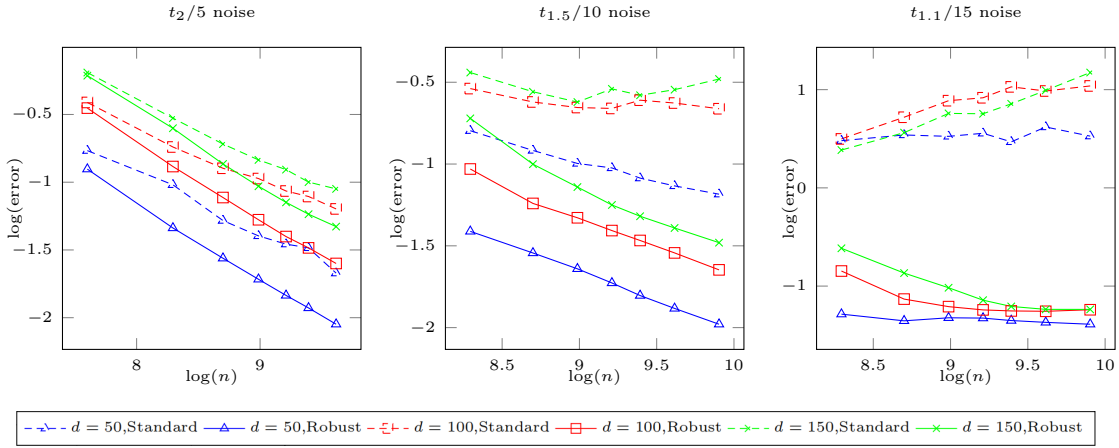

In matrix completion model, let $\Theta^{\star}=V_{5}V_{5}^T/\sqrt{5}$ where $V_{5}$ is top $5$ eigenvectors of $d$-dimensional sample covariance matrix from $100$ i.i.d. standard Gaussian random vectors. We use almost the same algorithm (the ADMM method proposed by Fang et al. (2015)\cite{Fang}) as that of Fan et al. (2021)\cite{shrinkage}. The main difference is that we adapt $\Theta_{i, j}^{n}=\sum_{i=1}^{n} d_{1} d_{2} 1_{\left\{X_{i}=\sqrt{d_{1} d_{2}} e_{i} e_{j}^{T}\right\}} $, $\Theta_{i, j}^{s}=\sqrt{d_{1} d_{2}} \sum_{i=1}^{n} y_{i} 1_{\left\{X_{i}=\sqrt{d_{1} d_{2}} e_{i} e_{j}^{T}\right\}}$ in their algorithm and $X_{i}$ is sampled from $\{\sqrt{d_{1}d_{2}}\cdot e_{j}e_{k}^T\}_{j\in[d_{1}], k\in[d_{2}]}$. For computational convenience, Fan et al. (2021)\cite{shrinkage} proposed a robust cross-validation procedure without adhering to the derived rates of $\tau$ and $\lambda$. In this work, to demonstrate the phase transition of the statistical rate, we select $C_{1}\left(\frac{n}
{(d_{1}\vee d_{2})\log(d_{1}+d_{2})}\right)^{\frac{1}{\alpha}}$ and $C_{2}\left(\frac{(d_{1}\vee d_{2})\log(d_{1}+d_{2})}{n}\right)^{\frac{\alpha-1}{\alpha}}$ as $\tau$ and $\lambda$ respectively. $C_{1}$ and $C_{2}$ are the fixed constants for each line in Figure \ref{fig1}. Consider the scaled Student's $t_{\nu}$ distribution with $\nu \in \{1.1, 1.5, 2\}$ as the error distribution and take $\alpha=\nu-0.01$ in the simulation. The numerical results are presented based on the mean of $200$ independent repetitions. In  Figure \ref{fig1}, the slopes of the fitted lines via the robust procedure become lower as $\alpha$ decreases, which is in keeping with Theorem \ref{theorem1}. Besides, when the tail of the noise distribution is heavier, the robust estimator performs better than the standard procedure in which the responses are not clipped.\par

Next, we select the following set of functions as the link functions $\{f_{i}(\cdot): i \in [9]\}$ to verify the behavior of the robust estimator in (\ref{eq2.5}) with respect to the sample size:
\begin{figure}[H]
\centering
\subfigure{
\begin{tikzpicture}[scale=0.5]
\begin{axis}[legend columns=1,legend pos=north west, height=9cm, width=9.3cm,
    axis lines = left,
    xlabel = $x$,
    ylabel = {$f(x)$},
    line width=0.4mm,
    title=(a)]
\addplot [
    domain=-50:50,
    samples=100,
    color=red,]
{4*x*cos(5*x)^2)};
\addlegendentry{$4x\cos^{2}(5x)$}
\addplot [
    domain=-50:50,
    samples=100,
    color=blue,
    ]
    {4*x*sin(5*x)^2};
\addlegendentry{$4x\sin^{2}(5x)$}
\addplot [
    domain=-50:50,
    samples=100,
    color=green,
    ]
    {-5*x/(2+sin(x))};
\addlegendentry{$-5x/(2+\sin(x))$}
\end{axis}
\end{tikzpicture}}
\subfigure{
\begin{tikzpicture}[scale=0.5]
\begin{axis}[legend columns=1,legend pos=north west, height=9cm, width=9.3cm,
    axis lines = left,
    xlabel = $x$,
    ylabel = {$f(x)$},
    line width=0.4mm,
    title=(b)]
\addplot [
    domain=-10:10,
    samples=100,
    color=red,]
{4*x+exp(x)/(1+exp(x))};
\addlegendentry{$4x+\frac{\exp(x)}{1+\exp(x)}$}
\addplot [
    domain=-10:10,
    samples=300,
    color=blue,
    ]
    {2*x+exp(-x^2/7)};
\addlegendentry{$2x+\exp(-x^2/7)$}
\addplot [
    domain=-10:10,
    samples=100,
    color=green,]
{-x+5*cos(8*x)};
\addlegendentry{$-x+5\cos(8x)$}
\end{axis}
\end{tikzpicture}}
\subfigure{
\begin{tikzpicture}[scale=0.5]
\begin{axis}[legend columns=1,legend pos=north east, height=9cm, width=9.3cm,
    axis lines = left,
    xlabel = $x$,
    ylabel = {$f(x)$},
    line width=0.4mm,
    title=(c)]
\addplot [
    domain=-30:30,
    samples=100,
    color=red,]
{x+4*sin(7*x)};
\addlegendentry{$x+4\sin(7x)$}
\addplot [
    domain=-30:30,
    samples=100,
    color=blue,
    ]
    {-x+cos(3*x^2/2)};
\addlegendentry{$-x+\cos(3x^2/2)$}
\addplot [
    domain=-30:30,
    samples=100,
    color=green,
    ]
    {-2*x+4*sin(x^2/2)};
\addlegendentry{$-2x+4\sin(x^2/2)$}
\end{axis}
\end{tikzpicture}}
\caption{$(a): f_{1}(x)=4x\cos^{2}(5x)$, $f_{2}(x)=4x\sin^{2}(5x)$, $f_{3}(x)=-5x/(2+\sin(x))$; $(b): f_{4}(x)=4x+\frac{\exp(x)}{1+\exp(x)}$, $f_{5}(x)=2x+\exp(-x^2/7)$, $f_{6}(x)=-x+5\cos(8x)$ and $(c): f_{7}(x)=x+4\sin(7x)$, $f_{8}(x)=-x+\cos(3x^2/2)$, $f_{9}(x)=-2x+4\sin(x^2/2).$}\label{fig4}
\end{figure}
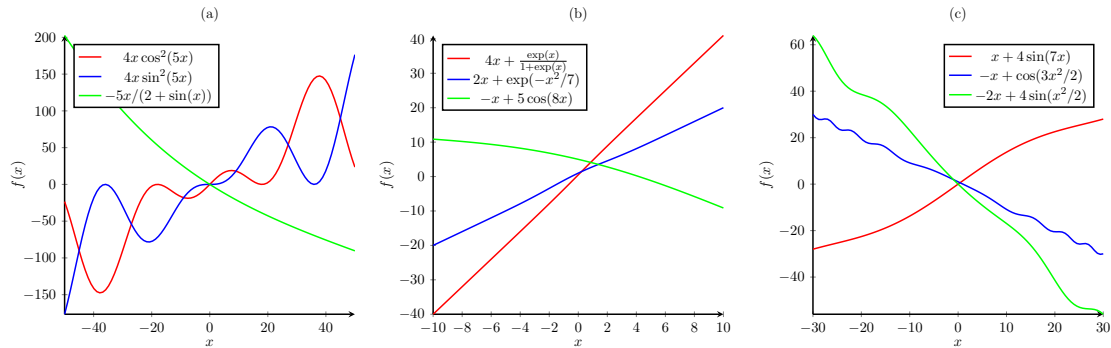
We set the dimensionality $d_{1}=200$ and for $\ell \in S_{k}$, $\left[\theta_{k}^{*}
\right]_{\ell}=\text { Uniform }(\{-1,1\})/\sqrt{s}$ where $S_{k}$ is the support of $\theta_{k}^{\star}$ chosen at random on $[d_{1}]$ with $\left|S_{k}\right|=s$. We use the distance
$$\rho\left(\widehat{\Theta}, \Theta^{*}\right)=\sqrt{\sum_{k=1}^{9}\min\left\{
\left\|\| \widehat{\theta}_{k}\|_{2}^{-1}\widehat{\theta}_{k}-\theta^{*}_{k}
\right\|_{2}^2,\left\|\| \widehat{\theta}_{k}\|_{2}^{-1}\widehat{\theta}
_{k}+\theta^{*}_{k}\right\|_{2}^2\right\}}$$
to measure the estimation error. Let the entries of $X$ and $\varepsilon$ i.i.d. follow $t_{5}$ distribution. $Z$ follows multivariate $t_{5}$ distribution where the precision matrix $\Omega$ is defined as $(\Omega)_{i,j}=0.5^{|i-j|}$. Inspired by Wang et al. (2021)\cite{tuning-free}, we solve the following adaptive equations to obtain truncation parameters $\left\{\tau_{j,k}^{(1)},\tau_{k,s}^{(2)}\right\}_{j \in [d_{1}]}^{k,s\in [d_{2}]}$ with computational efficiency:
$$
\sum_{i=1}^{n}\psi_{\tau_{j,k}^{(1)}}^{2}\left(y_{i}[S(X_{i})]_{j}z_{k}^{(i)}
\right)/\left(\tau_{j,k}^{(1)}\right)^{2}=10\log(d_{1}d_{2})\;\;\text{and} \;\; \sum_{i=1}^{n}\psi_{\tau_{k,s}^{(2)}}^{2}\left(z_{k}^{(i)}z_{s}^{(i)}
\right)/\left(\tau_{k,s}^{(2)}\right)^{2}=10\log(d_{2}).
$$
It is noteworthy that in the presence of heavy-tailed data and outliers, the above data-driven procedure can effectively select appropriate robustification parameters to truncate data. However, in Na et al. (2019)\cite{varying}, each truncation parameter needs to be adjusted by cross validation, which is inconvenient in practice. The numerical experiments are repeated 50 times.\par
\begin{table}[H]
    \centering
    \fontsize{8}{10}\selectfont
\caption{The logarithmic estimation error with respect to the sample size for $s=5$ and $10$.}
\begin{tabular}{|c|c|c|c|c|c|c|c|c|c|c|}
\hline
\diagbox{$s$}{$n$}& 10000 & 12500 & 15000 & 17500 & 20000 & 22500 & 25000 & 30000 & 35000 &  \\
\cline{1-11} \multirow{2}{*}{$5$} & $0.5647$ & $0.4662$ & $0.3829$ & $0.3095$ & $0.2620$ & $0.2100$ & $0.1540$ & $0.0849$ & $0.0046$ & Robust \\
\cline{2-11}  &0.7570 & 0.7019 & 0.6420 &  0.5759 & 0.5004 & 0.4340 & 0.4033 & 0.3709  & 0.2569 &Standard \\
\cline{1-11} \multirow{2}{*}{$10$} & 0.7144 & 0.6293 &0.5547&0.4779 &0.3928& 0.3458 & 0.2835& 0.1903 &0.1207 &Robust\\
\cline{2-11}  & 0.8801 & 0.7942 & 0.7473 & 0.7197 & 0.6428 & 0.5947 & 0.5393
& 0.4899 & 0.3800&Standard\\
\hline
\end{tabular}
\label{tab1}
\end{table}
For each $s$, we gather all the data points $\left(\log(\rho(\widehat{\Theta},\Theta^{*})),n\right)$ of the robust procedure to fit the linear regression relationship (i.e. $\log(\rho(\widehat{\Theta}
,\Theta^{*}))=\beta_{0}+\beta_{1}\log(n)$). The fitting results are that for $s=5$, $\beta_{1}=-0.4425$ with multiple $R^2=0.9993$ and for $s=10$, $\beta_{1}=-0.4825$ with multiple $R^2=0.9960$. Therefore, Table \ref{tab1} corroborates the result of Theorem \ref{theorem2} and shows that our proposed estimator has smaller statistical error than the standard procedure which means that the truncation techniques in (\ref{eq2.5}) and (\ref{eq1.6}) are not adopted.
\section{Concluding remarks}\label{sec5}
In this article, we extend Fan et al. (2021)\cite{shrinkage}'s work to the finite mean setting for heavy-tailed noise and observe a phase transition phenomenon by theory and experiment. Moreover, for high-dimensional varying index coefficient model, our proposed estimator is superior to Na et al. (2019)\cite{varying}'s robust estimator in two aspects. First, it allows the covariate and response to have bounded fourth moment. Second, tuning parameters via the data-driven procedure offers significant advantages in convenience and computing efficiency. The numerical experiments show that the improved estimator consistently performs better than the standard procedure and has consistency with the theoretical result. It is interesting that according to Sun et al. (2020)\cite{Adaptive}, the proposed estimator for part (a) of Theorem 2 in Fan et al. (2021)\cite{shrinkage} has also a tight phase transition phenomenon by following the proof of Theorem \ref{theorem1}. This strongly implies that the theoretical rate of Theorem \ref{theorem1} is sharp and we leave it as future research.

\section*{Acknowledgement}
We thank the editor and two anonymous reviewers for detailed and insight comments. This work was supported by grants from the NSF of China (Grant No.11731012), Ten Thousands Talents Plan of Zhejiang Province (Grant No. 2018R52042) and the Fundamental Research Funds for the Central Universities.

\section*{\textbf{Appendix A. Supplementary data}}
Some additional simulation results and the proofs of two theorems are contained in the supplementary material.

\end{document}